\documentclass[11pt]{article}
\usepackage{float}
\usepackage{fullpage}
\usepackage{calc}
\usepackage{amsmath}
\usepackage{tabularx}
\usepackage{amsfonts}
\usepackage{amsthm}
\usepackage{epsf}
\usepackage[dvips]{epsfig}
\usepackage{epsfig}
\begin{document}
\newcounter{claimcounter}
\newtheorem{theorem}{Theorem}[section]
\newtheorem{corollary}[theorem]{Corollary}
\newtheorem{proposition}[theorem]{Proposition}
\newtheorem{conjecture}[theorem]{Conjecture}
\newtheorem{lemma}[theorem]{Lemma}
\newtheorem{claim}[claimcounter]{Claim}
{\theoremstyle{definition} \newtheorem{definition}[theorem]{Definition}}
{\theoremstyle{definition} \newtheorem{example}[theorem]{Example}}
\newenvironment{notation}[1][Notation.]{\begin{trivlist}
\item[\hskip \labelsep {\bfseries #1}\\]}{\end{trivlist}}

\renewcommand{\thefigure}{\arabic{section}.\arabic{theorem}}
\newcommand{\stirlingtwo}[2]{\genfrac{\lbrace}{\rbrace}{0pt}{}{#1}{#2}}
\newcommand{\ceil}[1]{\left\lceil#1\right\rceil}
\newcommand{\floor}[1]{\left\lfloor#1\right\rfloor}
\newcommand{\vspan}{\operatorname{span}}
\newcommand{\rank}{\operatorname{rank}}
\newcommand{\img}{\operatorname{img}}
\newcommand{\rt}{\operatorname{root}}
\newcommand{\RT}{\operatorname{ROOT}}
\newcommand{\stcmp}{\operatorname{sc}}

\thispagestyle{empty}

\begin{center}
	{\LARGE State Complexity and the Monoid of Transformations\\ 
	of a Finite Set}\\
	\vspace*{.25in}
	\begin{tabular}{ccc}
		{\large Bryan Krawetz} &{\large John Lawrence} &{\large Jeffrey Shallit}\\
		{\em School of Computer Science} &{\em Department of Pure Mathematics} &{\em School of Computer Science}\\
		{\tt bakrawet@uwaterloo.ca} &{\tt jwlawren@math.uwaterloo.ca} &{\tt shallit@uwaterloo.ca}\\
	\end{tabular}\\
	\vspace*{.25in}
	{\em University of Waterloo\\ Waterloo, Ontario, Canada N2L 3G1}\\
\end{center}

\begin{abstract}
In this paper we consider the state complexity of an operation on formal languages, $\rt(L)$.
This naturally entails the discussion of the monoid of transformations of a finite set. 
We obtain good upper and lower bounds on the state complexity of $\rt(L)$ over alphabets
of all sizes.
\end{abstract}

\pagenumbering{arabic}

\section{Introduction}
\label{sec-intro}
A {\em deterministic finite automaton}, or {\em DFA}, is a 5-tuple ${\mathcal A}=(Q,\Sigma,\delta,q_0,F)$,
where $Q$ is a finite non-empty set of states, $\Sigma$ is the finite input alphabet, 
$\delta:Q\times\Sigma\rightarrow Q$ is the transition function, $q_0 \in Q$ is the initial state, 
and $F \subseteq Q$ is the set of final states. We assume that $\delta$ is defined on all elements of its
domain.  The domain of $\delta$ can be extended in the obvious way to $Q\times\Sigma^*$, where $\Sigma^*$ is
the free monoid over the alphabet $\Sigma$.  For a DFA ${\mathcal A}$, 
the set $L({\mathcal A}) = \{ w \in \Sigma^* \;:\; \delta(q_0,w) \in F\}$ is
said to be the language recognized by ${\mathcal A}$.\\

The {\em state complexity} of a regular language $L\subseteq\Sigma^*$, denoted $\stcmp(L)$, is defined 
as the size (the number of states) of the smallest DFA recognizing $L$. The state complexity of 
various operations on regular languages, such as union, concatenation, and Kleene closure, has been 
studied extensively; see, for example, \cite{Yu99,YZS94}.\\

In this paper we examine a less familiar operation, namely $\rt(L)$, which is given by
	\[\rt(L) = \{ w \in \Sigma^* : \exists n \geq 1 \text{ such that } w^n \in L\}.\]
Note that this operation is not the same as the $\RT(L)$ operation studied by Horv\'ath, Leupold, 
and Lischke \cite{HLL02}. The study of the $\rt(L)$ operation requires us to examine 
the connections between finite automata and algebra.\\

For a finite set $Q$, a function $f:Q \rightarrow Q$ is called a {\em transformation}. The
set of all transformations of $Q$ is denoted $Q^Q$. For transformations $f,g \in Q^Q$, 
their composition is written $fg$, and is given by $(fg)(q) = g(f(q))$, for all $q \in Q$.  
Together, the set $Q^Q$ and the composition operator form a monoid.\\

Transformations and their monoids have been studied in some detail
by D\'enes (whose work is summarized in \cite{Dene72}), and Salomaa \cite{Salo63,Salo03}.  
D\'enes investigates several algebraic and combinatorial properties of transformations, 
while much of Salomaa's work is concerned with subsets that generate the full 
monoid of transformations.\\

Let $L$ be a language and ${\mathcal A} = (Q,\Sigma,\delta,q_0,F)$ a DFA such that $L = L({\mathcal A})$.
For $w \in \Sigma^*$, define $\delta_w(q) = \delta(q,w)$, for all $q \in Q$. Then
$\delta_w$ is a transformation of $Q$. If we denote the empty word by $\epsilon$, 
then $\delta_{\epsilon}$ is the identity transformation.

\begin{theorem}
\label{thm-root-reg}
For a language $L$ and a DFA ${\mathcal A} = (Q,\Sigma,\delta,q_0,F)$ with $L = L({\mathcal A})$,
define the DFA ${\mathcal A}' = (Q^Q,\Sigma,\delta',q_0',F')$ where $q_0' = \delta_{\epsilon}$, 
$F' = \{f : \exists n\geq 1 \text{ such that } f^n(q_0) \in F\}$, and $\delta'$ is given by
\begin{align*}
	\delta'(f,a) = f\delta_a, \text{ for all $f\in Q^Q$ and $a\in\Sigma$}.
\end{align*}
Then $\rt(L) = L(\mathcal{A}')$.
\end{theorem}

\begin{proof}
An easy induction on $|w|$, $w \in \Sigma^*$, proves that $\delta'(q_0',w) = \delta_w$. 
Then
\begin{align*}
	x \in \rt(L) &\Leftrightarrow \exists n \geq 1\ :\ x^n \in L\\
		&\Leftrightarrow \delta_x(q_0) \in F'\\
		&\Leftrightarrow \delta'(q_0',x) \in F'.
\end{align*}
\end{proof}

In addition to giving us a construction for a DFA recognizing $\rt(L)$, the
above result shows that this operation preserves regularity, that is,
if $L$ is regular, then $\rt(L)$ is regular.  Zhang \cite{Zhan99} used a similar technique to 
characterize regularity-preserving operations. To recognize the image of a language 
under an operation, Zhang constructs a new automaton with states based on Boolean matrices.  
These matrices represent the transformations of states in the original automaton.\\

The result of Theorem \ref{thm-root-reg} also allows us to give our first 
bound on the state complexity of $\rt(L)$.

\begin{corollary}
\label{cor-root-bound}
For regular language $L$, if $\stcmp(L) = n$ then $\stcmp(\rt(L)) \leq n^n$.
\end{corollary}

\begin{proof}
This is immediate from the construction given in Theorem \ref{thm-root-reg}.
\end{proof}

In the remainder of this paper we improve on this upper bound, and give a non-trivial lower bound for
the worst-case blow-up of the state complexity of $\rt(L)$, for alphabets of all sizes. These upper
and lower bounds demonstrate that a simple, intuitive operation that preserves regularity
can increase the state complexity of a language from $n$ to nearly $n^n$, even over binary alphabets.\\

Our main results are given in Corollary \ref{cor-2-bound}, Corollary \ref{cor-2-prime-max}, 
and Theorem \ref{thm-3-max}.

\section{Unary languages}
\label{sec-unary}

In the case of unary regular languages, it turns out that the state 
complexity of the root of a language is bounded by the state complexity of
the original language.

\begin{proposition}
If $L$ is a unary regular language, then $\stcmp(\rt(L)) \leq \stcmp(L)$.
This bound is tight.
\end{proposition}

The idea of the following proof is that given a particular DFA recognizing $L$, we 
can modify it by adding states to the set of final states. The resulting DFA will 
recognize the language $\rt(L)$.

\begin{proof}
Let $\Sigma = \{a\}$ be the alphabet of $L$. Since $L$ is regular and unary, 
there exists a DFA $\mathcal{A}$ recognizing $L$, such that
$\mathcal{A} = (\{q_0,\ldots,q_{n-1}\},\{a\},\delta,q_0,F)$, where
\begin{align*}
	\delta(q_i,a) &= q_{i+1}, \text{for all $0 \leq i < n-1$},
\intertext{and} 
	\delta(q_{n-1},a) &= q_j, \text{for some $0 \leq j \leq n-1$}.
\end{align*}
We call the states $q_0,\ldots,q_{j-1}$ the {\em tail}, and the states $q_j,\ldots,q_{n-1}$ the {\em loop}.\\

Notice that $\rt(L) = \{a^s \in \Sigma^* : s\ |\ t, a^{t}\in L\}$. For all strings $a^t \in L$,
we have some $k \geq 0$ and some $b \leq n-1$ such that $t = kl + b$, where $l = n - j$ 
is the number of states in the loop.  Let $s = lm+c$ for some $m \geq 0$ and some $0 \leq c < l$. 
Then
\begin{align*}
	s\ |\ t\ &\Leftrightarrow\ \exists r\ :\ lk+b = r(lm+c)\\
		&\Leftrightarrow\ \exists r\ :\ lk - rlm = rc-b\\
		&\Leftrightarrow\ \exists r\ :\ \gcd(l, -lm)\ |\ rc-b &\text{(by Theorem 4.3.1 of \cite{BS96})}\\
		&\Leftrightarrow\ \exists r\ :\ l\ |\ rc-b\\
		&\Leftrightarrow\ \exists r,v\ :\ rc-b = lv\\
		&\Leftrightarrow\ \exists r,v\ :\ rc-lv = b\\
		&\Leftrightarrow\ \gcd(l,c)\ |\ b. &\text{(by Theorem 4.3.1 of \cite{BS96})}
\end{align*}
It follows that the set of divisors of the numbers of the form $kl + b$, $k\geq0$, $b \leq n-1$ 
is as follows:
	\[\{ lm +c \in {\mathbb Z} : m\geq0, \gcd(l,c)\ |\ b\}.\]
These divisors can be recognized by changing the corresponding states into final
states. Therefore, $\stcmp(\rt(L)) \leq \stcmp(L)$.\\

To show that this bound is tight, for $n \geq 2$ consider the language $L_n = \{a^{n-2}\}$. 
Under the Myhill-Nerode equivalence relation \cite{HU79}, no two strings in the set 
$\{\epsilon,a,a^2,\ldots,a^{n-1}\}$ are equivalent.  All other strings in $\Sigma^*$ are 
equivalent to $a^{n-1}$. This gives $\stcmp(L_n) = n$.  Furthermore, since $a^{n-2}$ is the 
longest word in $\rt(L_n)$, $\delta(q_0,a^{n-2})$ cannot be a state in the loop.  It follows 
that we require exactly $n-1$ states in the tail plus a single, non-final state in the loop. 
Hence $\stcmp(\rt(L_n)) = n$.  Therefore the bound is tight.
\end{proof}

\section{Languages on larger alphabets}
\label{sec-larger}
For a regular language $L \subseteq \Sigma^*$, if ${\mathcal A}$ is the minimal DFA such that $L = L({\mathcal A})$,
then as we saw in Section \ref{sec-intro}, based on the set of all transformations of the states of ${\mathcal A}$, 
we can construct an automaton ${\mathcal A}'$ as in Theorem \ref{thm-root-reg}, to recognize $\rt(L)$. 
Though this new DFA, $\mathcal{A}'$, has all transformations of $Q$ as its states, it is easy to see that the 
only reachable states are those that are a composition of the transformations 
$\delta_{a_1},\ldots,\delta_{a_m}$, where $a_1,\ldots,a_m \in \Sigma$.  This set of elements, 
$\delta_{a_1},\ldots,\delta_{a_m}$, and all of their compositions form the transformation monoid of 
$\mathcal{A}$.  We use this fact to improve on the upper bound of the state complexity of $\rt(L)$.

\begin{corollary}
\label{cor-root-monoid-bound}
For a regular language $L$, let $\mathcal{A}$ be the smallest DFA recognizing $L$.  Then if $M$ is 
the transformation monoid of $\mathcal{A}$, we have that $\stcmp(\rt(L)) \leq |M|$.
\end{corollary}

\begin{proof}
In Theorem \ref{thm-root-reg}, the only reachable states in the construction of $\mathcal{A}'$ 
are those that belong to the transformation monoid of $\mathcal{A}$.
\end{proof}

Define $Z_n = \{1,2,\ldots,n\}$. Now define $T_n = Z_n^{Z_n}$, the set of 
transformations of $Z_n$, and $S_n \subseteq T_n$ as the set of permutations
of $Z_n$. For $\gamma \in T_n$ we write 
\[ \gamma = \begin{pmatrix} 1 &2 &\cdots &n\\ \gamma(1) &\gamma(2) &\cdots &\gamma(n)
	\end{pmatrix}. \]

\begin{definition}
If $M \subseteq T_n$ is the set of all compositions of the transformations $f_1,\ldots,f_m \in T_n$, 
then we say that $\{f_1,\ldots,f_m\}$ {\em generates} $M$.
\end{definition}

\begin{definition}
For $\gamma \in T_n$, define the {\em image} of $\gamma$ by
$\img(\gamma) = \{ y \in Z_n : y = \gamma(z) \text{ for some } z \in Z_n\}$.
\end{definition}

\begin{definition}
For $\gamma \in T_n$, define the {\em rank} of $\gamma$ as the number of distinct elements 
in the image of $\gamma$, and denote it by $\rank(\gamma)$.
\end{definition}

The relationship between the state complexity of a language and the transformation monoid naturally leads 
to the question of how large a submonoid of $T_n$ can be generated by $m$ elements, where
$m$ is a positive integer.  In connection with the study of Landau's function (for a survey see 
\cite{Mill87,Nico97}), Szalay \cite{Szal80} showed that, for $m=1$, 
the largest submonoid of $T_n$ has size
	\[ \exp\left\{\sqrt{n\left(\log n + \log\log n -1 + \frac{\log\log n -2 + o(1)}{\log n}\right)}\right\}.\] 

In the case where $m \geq 3$, the results are well known.

\begin{lemma}
\label{lem-semi-gen}
For $n \geq 3$, suppose $H \subseteq T_n$ such that $H$ generates  $T_n$. 
Then $|H| \geq 3$.  Furthermore, $|H| = 3$ if and only if $H$ can be written as 
$H = \{\alpha,\beta,\gamma\}$, where $\{\alpha, \beta\}$ generates $S_n$ and 
$\rank(\gamma) = n-1$.
\end{lemma}

For a proof of this lemma, see D\'enes \cite{Dene66}. This gives us the result 
that the largest submonoid generated by three elements has the full size $n^n$.\\

Contrary to the case for $m=1$ and $m\geq3$, it seems that only recently there has been 
any interest in determining the largest submonoid on two generators.  
Significant progress has been made in this area by Holzer and K\"onig \cite{HK02,HK03}, 
and, independently, by Krawetz, Lawrence, and Shallit \cite{KLS03}. The
results of Holzer and K\"onig are summarized here.\\ 

For coprime integers $k,l \geq 2$, where $k+l=n$, let $\alpha = (1\ 2\ \cdots\ k)(k+1\ \ k+2\ \ \cdots\ n)$
be a permutation of $Z_n$ composed of two cycles, one of length $k$, the other
of length $l$. Define $U_{k,l}$ to be the set of all transformations $\gamma \in T_n$ 
where exactly one of the following is true:
\begin{enumerate}
	\item $\gamma = \alpha^m$ for some positive integer $m$;
	\item For some $i \in \{1,\cdots, k\}$ and some $j \in \{k+1,\cdots, n\}$ we
		have that $\gamma(i) = \gamma(k)$ and for some $m \in \{k+1,\cdots, n\}$ we
		have that $m \not\in \img(\gamma)$.
\end{enumerate}

Let $\pi_1 = (1\ 2\ \cdots\ k)$ be an element of $S_{n-1}$, and let $\pi_2 \in S_{n-1}$ be
a permutation such that $\pi_1$ and $\pi_2$ generate $S_{n-1}$.  Now define $\beta \in T_n$ by
\begin{align*}
\beta &= \begin{pmatrix}
		1 &2 &\cdots &n-1 &n\\
		\pi_2(1) &\pi_2(2) &\cdots &\pi_2(n-1) &\pi_2(1)
	\end{pmatrix}.
\end{align*}

\begin{lemma}[Holzer and K\"onig]
\label{lem-u-monoid}
The set $U_{k,l}$ is a submonoid of $T_n$ and is generated by 
$\{\alpha,\beta\}$.
\end{lemma}

It is worth noting that in their definition of $U_{k,l}$,  Holzer and K\"onig allow $k=1$ and $l=1$.
They show implicitly, however, that the size of the monoid in these degenerate cases is too
small to be of any consequence here.

\begin{theorem}[Holzer and K\"onig]
\label{thm-monoid-size}
For $n \geq 7$, there exist coprime integers $k$,$l$ such that $n = k + l$ and 
	\[|U_{k,l}| \geq n^n\left(1 - \sqrt{2}\left(\frac{2}{e}\right)^{\frac{n}{2}}e^{\frac{1}{12}} 
		- \sqrt{8}\frac{1}{\sqrt{n}}e^{\frac{1}{12}}\right).\]
\end{theorem}

In addition to a lower bound on the size of the largest two-generated monoid, Theorem \ref{thm-monoid-size} 
gives us the existence of a sequence of two-generated monoids whose size approaches $n^n$ as $n$ 
tends toward infinity.  Similar results were obtained independently by Krawetz, Lawrence, 
and Shallit \cite{KLS03}.\\

More recently, Holzer and K\"onig \cite{HK03} proved the following result regarding
the maximality of monoids of the form $U_{k,l}$. 

\begin{theorem}[Holzer and K\"onig]
\label{thm-max-prime}
For all prime numbers $n \geq 7$, there exist coprime integers $k$,$l \geq 2$ such that $k+l=n$ 
and $U_{k,l}$ is the largest two-generated submonoid of $T_n$.
\end{theorem}

They also stated the following conjecture.

\begin{conjecture}[Holzer and K\"onig]
\label{con-ukl-max}
For any $n \geq 7$, there exist coprime integers $k$,$l \geq 2$ such that $k+l=n$ and $U_{k,l}$ is the 
largest two-generated submonoid of $T_n$.
\end{conjecture} 

Since the connection between the state complexity of $\rt(L)$ and the transformation monoid of $L$ has
been established in Corollary \ref{cor-root-monoid-bound}, we can take advantage of the results of 
Theorem \ref{thm-monoid-size} and Theorem \ref{thm-max-prime} if we can construct a language based on a monoid.
By associating an alphabet with the generators of a monoid, we can define a transition 
function for a DFA. The definition of the DFA is then completed by choosing a start 
state and a set of final states.  This construction is given more formally below.\\
 
Let $n,m$ be integers with $n,m\geq1$. For a set of transformations  $X=\{\alpha_1,\ldots,\alpha_m\}$,
let $M \subseteq T_n$ denote the monoid generated by $X$.  Then a {\em DFA based on $X$} is a 
DFA ${\mathcal M} = (Z_n,\Sigma,\delta,z_0,F)$, where $|\Sigma| \geq m$, $z_0 \in Z_n$, 
$F \subseteq Z_n$, and $\delta$ is given by 
	\[\delta_{a} = \Psi(a),\ \text{for all $a \in \Sigma$,}\]
for some map $\Psi:\Sigma \rightarrow X \cup \{\delta_{\epsilon}\}$ that is surjective on $X$.

\begin{proposition}
Let ${\mathcal M} = (Z_n,\Sigma,\delta,z_0,F)$ be a DFA. Then $M$ is the transformation monoid of ${\mathcal M}$
if and only if ${\mathcal M}$ is based $X$, for some $X\subseteq T_n$ that generates the monoid $M$. 
\end{proposition}

\begin{proof}
For a DFA ${\mathcal M}$ based on $X$, the fact that $M$ is the transformation monoid of ${\mathcal M}$
is immediate from the construction.  For any DFA ${\mathcal M}$ that has $M$ as its transformation 
monoid, we have that the set $\{ \delta_a \in T_n : a \in \Sigma\}$ generates $M$.  Then we can simply
take $\Psi$ given by $\Psi(a) = \delta_a$, for all $a \in \Sigma$.
\end{proof}

In particular, let ${\mathcal A}_{\Psi,X} = (Z_n,\Sigma,\delta,z_0,F)$ denote the DFA based on $X$
when $z_0 = 1$, $F = \{1\}$, and $\Psi$ is bijective on an $m$-element subset of $\Sigma$, with  
all other elements of $\Sigma$ mapped to $\delta_{\epsilon}$.  If $\Psi_1$ and $\Psi_2$ are
maps over the same domain, then ${\mathcal A}_{\Psi_1,X}$ is isomorphic to ${\mathcal A}_{\Psi_2,X}$, 
up to a renaming of the states and alphabet symbols.  For this reason, we will often denote this 
DFA simply by ${\mathcal A}_{\Sigma,X}$.\\

\begin{example} 
Let $Y = \{\alpha,\beta\}$, where 
\[\alpha = \begin{pmatrix} 1 &2 &3 &4 &5\\2 &1 &4 &5 &3\end{pmatrix},\ \text{and}\
\beta = \begin{pmatrix} 1 &2 &3 &4 &5\\2 &3 &4 &1 &2\end{pmatrix}.\]
Define $\Phi$ by $\Phi(a) = \alpha$ and $\Phi(b) = \beta$.  Then Figure \ref{fig-dfa} 
depicts the DFA ${\mathcal A}_{\Phi,Y}$.\\

\begin{figure}[H]
	\addtocounter{theorem}{1}
	\begin{center}
		\epsfig{file=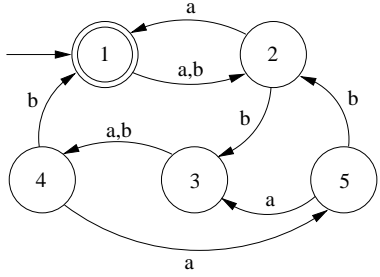,scale=0.75}
		\caption{The automaton $\mathcal{A}_{\Phi,Y}$.}
		\label{fig-dfa}
	\end{center}
\end{figure}
\end{example}

For $n \geq 5$, define $X_n \subseteq T_n$ to be a subset of transformations, and let $M_n$ denote the
monoid generated by $X_n$. In particular, define $X_{k,l} = \{\alpha,\beta\}$, where $\alpha$ and $\beta$ are
as in Lemma \ref{lem-u-monoid}. Then $X_{k,l}$ generates $U_{k,l}$.  We now state our main result 
concerning the state complexity of $\rt(L({\mathcal A}_{\Sigma,X_n}))$.

\begin{theorem}
\label{thm-min-dfa}
If $U_{k,l} \subseteq M_n$, for some coprime integers $k\geq2$ and $l\geq3$, with $k+l=n$, then the 
minimal DFA recognizing $\rt(L({\mathcal A}_{\Sigma,X_n}))$ has $|M_n| - \binom{n}{2}$ states.
\end{theorem}

Before we are ready to prove this theorem, we must state a few more definitions and lemmas.

\begin{definition}
Let $\rho \in T_n$. For any $i,j,k$, if $\rho(i) = k = \rho(j)$ implies that $i = j$, then
we say that $k$ is {\em unique}.
\end{definition}

\begin{definition}
Let $\rho \in T_n$ have rank $2$, with $\img(\rho) = \{i,j\}$. Then by the {\em complement} 
of $\rho$, we mean the transformation $\overline{\rho} \in T_n$, where 
\begin{align*}
	\overline{\rho}(k) = \begin{cases}
		i,& \text{if }\rho(k) = j;\\
		j,& \text{if }\rho(k) = i.
	\end{cases}
\end{align*}
\end{definition}

For example, if 
$\rho = \begin{pmatrix}
		1 &2 &3 &\cdots &n-1 &n\\
		3 &3 &2 &\cdots &2   &2
	\end{pmatrix}$, then
$\overline{\rho} = \begin{pmatrix}
		1 &2 &3 &\cdots &n-1 &n\\
		2 &2 &3 &\cdots &3   &3
	\end{pmatrix}$.\\

It is easy to see that, in general, $\rho$ and $\overline{\rho}$ have the same rank, 
and $\overline{\overline{\rho}} = \rho$.\\

For an automaton ${\mathcal M} = (Z_n,\Sigma,\delta,z_0,F)$, define the DFA 
${\mathcal M}^* = (M_n,\Sigma,\delta',\delta_{\epsilon},F')$, where $\delta_{\epsilon}$ 
is the identity element of $T_n$, $\delta'(\eta,a) = \eta\delta_a$ 
for all $\eta \in M_n$, for all $a \in \Sigma$, and $F' = \{\eta \in T_n : \eta(z_0) = z_0 \}$.  
Then $L({\mathcal M}^*) =$ $\rt(L({\mathcal M}))$.\\

For $\eta$, $\theta \in M_n$, with $\eta \not= \theta$, note that $\eta$ and $\theta$ 
are equivalent states if and only if for all $\rho \in \Sigma^*$
we have that 
	\[\delta'(\eta,\rho) \in F' \Leftrightarrow \delta'(\theta,\rho) \in F'.\]  
However, since $\delta'(\eta,\rho) = \eta\delta_{\rho}$,
this is equivalent to saying that $\eta$ and $\theta$ are equivalent states in $M$ if and only if 
for all $\rho \in M_n$, we have
	\[\eta\rho \in F' \Leftrightarrow \theta\rho \in F'.\]

\begin{lemma}
\label{lem-equiv-1}
Let $Y \subseteq T_n$ generate $M_n$, and let ${\mathcal M} = (Z_n,\Sigma,\delta,z_0,F)$ be an 
automaton based on $Y$ such that $z_0\in F$.  Let $\eta$, $\theta \in M_n$, 
with $\eta \not= \theta$ and $\rank(\eta) = 2$.  
If $\eta(z_0)$ is unique in the image of $\eta$, and $\overline{\eta} = \theta$, 
then $\eta$ and $\theta$ are equivalent states in ${\mathcal M}^*$.
\end{lemma}

\begin{proof}
We have
\begin{align*}
	\eta &= \begin{pmatrix}
		1 &2 &\cdots &z_0-1 &z_0 &z_0+1 &\cdots &n\\
		j &j &\cdots &j     &i   &j     &\cdots &j
	\end{pmatrix},
	\intertext{and} 
	\theta &= \begin{pmatrix}
		1 &2 &\cdots &z_0-1 &z_0 &z_0+1 &\cdots &n\\
		i &i &\cdots &i     &j   &i     &\cdots &i
	\end{pmatrix},
\end{align*}
for some $i \not= j$.\\

If $\eta(z_0) \in F$, then $\eta \in F'$.  If $\theta(z_0) \in F$, then $\theta \in F'$.
Otherwise, $\theta(z_0) \in Z_n\backslash\{z_0\}$.
Since $\theta(z) = \eta(z_0) \in F$ for all $z \in Z_n\backslash\{z_0\}$, we have that 
$\theta^2(z_0) \in F$, and hence $\theta \in F'$.  Similarly, $\theta(z_0) \in F$ 
implies that $\eta,\theta \in F'$.  Furthermore, if $\img(\eta) \cap F = \emptyset$, then 
$\img(\eta^n) \cap F = \emptyset$, for all $n$, and hence $\eta^n(z_0) \not\in F$ so that 
$\eta \not\in F'$.  Since $\img(\eta) = \img(\theta)$, this gives $\theta \not\in F'$.
Therefore $\eta \in F'$ if and only if $\theta \in F'$.\\

Let $\rho \in M_n$. Since $\eta$ and $\theta$ have rank $2$, we must have that 
$\eta\rho$ and $\theta\rho$ have rank $\leq 2$.
If $\eta\rho$ has rank $2$, then $\rho(i) \not= \rho(j)$, so that $\overline{\eta\rho} = \theta\rho$. 
Hence $\eta\rho \in F'$ if and only if $\theta\rho \in F'$. 
The argument is the same for the case where $\theta\rho$ has rank $2$.
Now, if $\eta\rho$ has rank $1$, then we must have $\rho = \rho'\sigma\rho''$, where
$\sigma(s) = \sigma(t)$ for some $s$ and $t$, and where $\rho'$ is a permutation such that either 
$\eta\rho'(z_0) = s$ and $\eta\rho'(z) = t$, for all $z \not= z_0$, or 
$\eta\rho'(z_0) = t$ and $\eta\rho'(z) = s$, for all $z \not= z_0$.  Without loss of generality, assume the former.
Then clearly $\theta\rho'(z_0) = t$ and $\theta\rho'(z) = s$, for all $z \not= z_0$. It follows
that $\eta\rho'\sigma = \theta\rho'\sigma$, so that $\eta\rho = \theta\rho$.\\

Therefore $\eta$ and $\theta$ are equivalent states.
\end{proof}

\begin{lemma}
\label{lem-equiv-2}
Let $Y \subseteq T_n$ generate $M_n$, and let ${\mathcal M} = (Z_n,\Sigma,\delta,z_0,F)$ be an 
automaton based on $Y$, such that $z_0\not\in F$.  Let $\eta$, $\theta \in M_n$, with $\eta(z_0)$ unique in the image of 
$\eta$, $\rank(\eta) = 2$, and $\img(\eta) = \img(\theta)$.  If $\theta(z_0) = \eta(z_0)$,
and if $\theta(z) = \eta(z_0)$ implies that $z \in F$, then $\eta$ and $\theta$ are 
equivalent states in ${\mathcal M}^*$.
\end{lemma}

\begin{proof}
If $\eta(z_0) \in F$, then $\eta \in F'$. Since $\theta(z_0) = \eta(z_0)$, it follows that
$\theta \in F'$. Now suppose that $\eta(z_0) \not\in F$. If $\eta \in F'$, then since $\rank(\eta) = 2$, we must 
have that $\eta^2(z_0) \in F$. Now $\theta(z_0) \not\in F$, so $\theta(\theta(z_0)) \not= \eta(z_0)$.
But since $\rank(\theta) = 2$, and $\img(\theta) = \img(\eta)$, it follows that $\theta^2(z_0) = \eta^2(z_0)$.
Hence $\theta \in F'$.  If $\eta \not\in F'$, then $\eta(z_0) = z_0$ or $\img(\eta) \cap F = \emptyset$.
In either case, this implies that $\theta \not\in F'$. Therefore $\eta \in F'$ if and only if $\theta \in F'$.\\  

Let $\rho \in M_n$. Then, following an argument similar to the one used in the proof of Lemma
\ref{lem-equiv-1}, we have that $\eta\rho \in F'$ if and only if $\theta\rho \in F'$. 
Therefore $\eta$ and $\theta$ are equivalent states.
\end{proof}

Now that we have a characterization of equivalent states in the general case for ${\mathcal M}^*$,
we turn our attention toward the specific case, for ${\mathcal A}_{\Sigma,X_n}^*$.

\begin{lemma}
\label{lem-not-equiv-1}
Let $\eta$, $\theta \in M_n$, with $\eta \not= \theta$. If $\rank(\eta) = 1$, 
then $\eta$ and $\theta$ are not equivalent states in ${\mathcal A}_{\Sigma,X_n}^*$.
\end{lemma}

\begin{proof}
Since $\eta$ has rank $1$, we have that $\img(\eta) = \{z_1\}$ for some $z_1$.
If $\eta(1) \not= \theta(1)$, then take $\rho \in U_{k,l}$ such that $\rho(z_1) = 2$, and $\rho(z) = 1$, 
for all $z \not= z_1$. Then $\img(\eta\rho) = \{2\}$, so that $\eta\rho \not\in F'$.  But 
$\theta\rho(1) = 1$, so that $\theta\rho \in F'$. Hence $\eta$ and $\theta$ are not equivalent.
If $\eta(1) = \theta(1)$, then $\rank(\theta)\not=1$ so that for some $z_2 \not= 1$ we have
$\theta(z_2) \not= z_1$. Take $\rho \in U_{k,l}$ such that $\rho(\theta(z_2)) = 1$, and $\rho(z) = z_2$, 
for all $z \not= \theta(z_2)$. Then $\img(\eta\rho) = \{z_2\}$, so that $\eta\rho \not\in F'$.  But 
$(\theta\rho)^2(1) = 1$, so that $\theta\rho \in F'$. Hence $\eta$ and $\theta$ are not equivalent.
\end{proof}

\begin{lemma}
\label{lem-not-equiv-2}
For $\eta$, $\theta \in M_n$, with $\eta \not= \theta$, let $\eta$ have rank $2$.  
Then $\eta$ and $\theta$ are equivalent states in ${\mathcal A}_{\Sigma,X_n}^*$ if and only if $\eta(1)$ is unique 
in the image of $\eta$, and $\overline{\eta} = \theta$.
\end{lemma}

\begin{proof}
Let $\img(\eta) = \{i,j\}$ for some $i,j$.  Without loss of generality, assume that $\eta(1) = i$.\\

Since $1 \in F$, Lemma \ref{lem-equiv-1} applies and gives the result in the forward direction. For the other
direction, we have two cases.\\

Case 1. $i$ is not unique in the image of $\eta$.
\begin{quote}
	Case 1.a. $i = \eta(1) \not= \theta(1)$.
	\begin{quote}
		Since $i$ is not unique, $\eta(z) = i$, for some $z \not = 1$. Choose $\rho \in U_{k,l}$ 
		such that $\rho(i) = z$, and $\rho(\theta(1)) = 1$.  Then $\eta\rho(1) = \eta\rho(z) = z$
		so that $(\eta\rho)^n(1) = z$, for all $n \geq 0$. This gives us $\eta\rho \not \in F'$.
		But $\theta\rho(1) = 1$, so that $\theta\rho \in F'$. 
	\end{quote}
	Case 1.b. $i = \eta(1) = \theta(1)$.
	\begin{quote}
		Case 1.b.i. $\img(\eta) = \img(\theta)$.
		\begin{quote}
			Since $\img(\eta) = \img(\theta)$ and $\eta \not = \theta$,
			then for some $z \not = 1$, we must have either $i = \eta(z) \not = \theta(z) = j$, or $j = \eta(z) 
			\not= \theta(z) = i$.  Without loss of generality, assume the former. 
			Then choose $\rho \in U_{k,l}$ such that $\rho(i) = z$, and $\rho(j) = 1$.  
			Then $(\eta\rho)^n(1) = z$, for all $n \geq 0$, and $(\theta\rho)^2(1) = 1$.
			This gives us $\eta\rho \not \in F'$ and $\theta\rho \in F'$. 
		\end{quote}

		Case 1.b.ii. $\img(\eta) \not= \img(\theta)$.
		\begin{quote}
			If $\rank(\theta) = 1$ then by  Lemma \ref{lem-not-equiv-1}, $\eta$ and $\theta$ 
			are not equivalent. Otherwise there exists some $z_1 \in Z_n$ such that 
			$\theta(z_1) \not\in\img(\eta)$.  Take $\rho \in U_{k,l}$ such that $\rho(\theta(z_1)) = 1$, 
			and $\rho(z) = 2$, for all $z \not= \theta(z_1)$.  Then $\eta\rho \not= \theta\rho$ 
			and $\rank(\eta\rho) = 1$, and so by Lemma \ref{lem-not-equiv-1}, $\eta\rho$ and $\theta\rho$ 
			are not equivalent. Hence $\eta$ and $\theta$ are not equivalent.
%
		\end{quote}
	\end{quote}
\end{quote}

Case 2. $i$ is unique in the image of $\eta$, and $\overline{\eta} \not= \theta$.
\begin{quote}
	If $\rank(\theta) = 1$ then by  Lemma \ref{lem-not-equiv-1}, $\eta$ and $\theta$ are not equivalent.
	If $\img(\eta) = \img(\theta)$, then since $\overline{\eta} \not= \theta$, we have that $\theta(1)$ is 
	not unique.  So we can reverse the roles of $\eta$ and $\theta$ and apply case $1$ to get the desired result.
	Assume then, that $\img(\eta) \not= \img(\theta)$, and $\rank(\theta) \geq 2$. Then
	the fact that $\eta$ and $\theta$ are not equivalent follows just as in case 1.b.ii.
\end{quote}

Therefore, $\eta$, and $\theta$ are equivalent states if and only if $\eta(1)$ is unique in 
the image of $\eta$, and $\overline{\eta} = \theta$.
\end{proof}

\begin{lemma}
\label{lem-not-equiv-3}
Let $\eta$, $\theta \in M_n$, with $\eta \not= \theta$, if $\eta$, $\theta$ have rank $\geq 3$, 
then $\eta$ and $\theta$ are not equivalent states in ${\mathcal A}_{\Sigma,X_n}^*$.
\end{lemma}

\begin{proof}
Since $\eta \not= \theta$, there exists some $z_1 \in Z_n$ such that 
$\eta(z_1) \not= \theta(z_1)$. Let $z_2 = \eta(z_1)$. Take $\rho \in U_{k,l}$ such that $\rho(z_2) = 1$, 
and $\rho(z) = 2$, for all $z \not= z_2$.  Since $\rank(\eta) \geq 3$, we have $rank(\eta\rho) = 2$.  If $\eta\rho(1)$
is not unique, then by Lemma \ref{lem-not-equiv-2}, $\eta\rho$ and $\theta\rho$ are not equivalent.
Hence $\eta$ and $\theta$ are not equivalent. If $\eta\rho(1)$ is unique, then it must be that $z_1 = 1$.  
Furthermore, since $\rank(\theta) \geq 3$, we cannot have $\theta(z) = z_2$ for all $z \not= 1$,
so we cannot have $\theta\rho(z) = 1$ for all $z \not= 1$. Therefore $\theta\rho \not= \overline{\eta\rho}$. 
Then by Lemma \ref{lem-not-equiv-2}, $\eta\rho$ and $\theta\rho$ are not equivalent. 
Hence $\eta$ and $\theta$ are not equivalent.
\end{proof}

We are now ready to prove Theorem \ref{thm-min-dfa}.

\begin{proof}[Proof (Theorem \ref{thm-min-dfa}).]
Lemma \ref{lem-not-equiv-1} -- \ref{lem-not-equiv-3} cover all possible cases for 
$\eta, \theta \in M_n$, $\eta \not= \theta$.  Therefore, two states are equivalent if and only if they
satisfy the hypothesis of Lemma \ref{lem-equiv-1}. There are $\binom{n}{2}$ such equivalence
classes in $U_{k,l} \subseteq M_n$, each containing exactly 2 elements. All other elements of $M_n$ are 
in equivalence classes by themselves. It follows that the minimal DFA recognizing $\rt(L({\mathcal A}_{\Sigma,X_n}))$ 
has $|M_n| - \binom{n}{2}$ states.
\end{proof}

Now that we have established a close relationship between $\stcmp(\rt(L))$ and the transformation monoid of the
the minimal automaton recognizing $L$, we can take advantage of results concerning the size of the largest 
monoids to give bounds on the worst-case blow-up of the state complexity of $\rt(L)$.  The following corollary 
gives a lower bound for alphabets of size two.  It also proves the existence of a sequence of regular binary languages
with state complexity $n$ whose root has a state complexity that approaches $n^n$ as $n$ increases without bound.\\

We now state our first main result.

\begin{corollary}
\label{cor-2-bound}
For $n\geq7$, there exists a regular language $L$ over an alphabet of size 2, with $\stcmp(L) \leq n$, such 
that 
\begin{align*}
	\stcmp(\rt(L)) \geq n^n\left(1 - \sqrt{2}\left(\frac{2}{e}\right)^{\frac{n}{2}}e^{\frac{1}{12}} 
		- \sqrt{8}\frac{1}{\sqrt{n}}e^{\frac{1}{12}}\right)-\binom{n}{2}.
\end{align*}
\end{corollary}

\begin{proof}

The result follows from a combination of Theorem \ref{thm-monoid-size} and Theorem \ref{thm-min-dfa}.
\end{proof}

Our results from Theorem \ref{thm-min-dfa} do not apply when $l = 2$. Unfortunately, 
Theorem \ref{thm-max-prime} does not exclude this possibility.  To guarantee that this fact
is of no consequence, we must show that not only is the monoid $U_{n-2,2}$ never the
largest, but that it is at least $\binom{n}{2}$ smaller than the largest monoid.  The
following lemma deals with this.

\begin{lemma}
\label{lem-max-by-gap}
For $n \geq 7$, we have that
\begin{align*}
	|U_{2,n-2}| - |U_{n-2,2}| &\geq \binom{n}{2}.
\end{align*}
\end{lemma}

Due to space constraints, the proof of this lemma has been relegated to the appendix.\\

The choice of start and final states in the construction of the DFA ${\mathcal A}_{\Sigma,X_n}$
is the best possible. The following theorem will show that for any other DFA with the 
same transition function, a different assignment of start and final states will not 
increase the state complexity of the language it recognizes.\\

\begin{theorem}
\label{thm-max-dfa}
Let $Y \subseteq T_n$ generate $M_n$, and let ${\mathcal M} = (Z_n,\Sigma,\delta,z_0,G)$ be an 
automaton based on $Y$. Then $\stcmp(\rt(L({\mathcal M}))) 
\leq \stcmp(\rt(L({\mathcal A}_{\Sigma,X_n}))$.
\end{theorem}

\begin{proof}
If $z_0 \in G$, then Lemma \ref{lem-equiv-1} applies.
It follows that there are at least $\binom{n}{2}$ pairs of equivalent states in $M^*$.
If $z_0 \not\in G$, then Lemma \ref{lem-equiv-2} applies, and again we have at least 
$\binom{n}{2}$ pairs of equivalent states in $M^*$. In either case, this gives
\[\stcmp(\rt(L({\mathcal M}))) \leq |M_n|-\binom{n}{2} \leq \stcmp(\rt(L({\mathcal A}_{\Sigma,X_n})).\]
\end{proof}

We now state our second main result.

\begin{corollary}
\label{cor-2-prime-max}
For prime numbers $n\geq7$, there exist positive, coprime integers $k\geq2$, $l\geq3$, with $k+l=n$, 
such that if $L$ is a language over an alphabet of size 2, with $\stcmp(L) \leq n$, then
$\stcmp(\rt(L)) \leq |U_{k,l}| - \binom{n}{2}$.  Furthermore, this bound is tight.
\end{corollary}

\begin{proof}
Let $U'$ denote the largest two-generated submonoid of $T_n$.  Then by Theorem \ref{thm-max-prime} 
and Lemma \ref{lem-max-by-gap}, we have that $U' = U_{k',l'}$ for some coprime integers 
$k'\geq2$, $l'\geq3$ with $k'+l' =n$.\\

Let ${\mathcal M}$ be the smallest DFA recognizing $L$, and let $M$ be the transformation monoid of 
${\mathcal M}$.  If $M$ is of the form $U_{k,l}$, with $k\geq2$, $l\geq3$, then $|U_{k,l}| \leq |U'|$. It follows
from Theorem \ref{thm-max-dfa} that $\stcmp(\rt(L)) \leq |U_{k,l}|-\binom{n}{2} \leq |U'|-\binom{n}{2}$.
If $M$ is of the form $U_{k,l}$, with $k=n-2$, $l=2$, then by Corollary \ref{cor-root-monoid-bound} and
Lemma \ref{lem-max-by-gap} we have $\stcmp(\rt(L)) \leq |U_{n-2,2}| \leq |U_{2,n-2}| -\binom{n}{2} \leq |U'|-\binom{n}{2}$
.\\

Let $V$ denote the largest two-generated submonoid 
of $T_n$ that is not of the form $U_{k,l}$ for some coprime integers $k\geq2$, $l\geq3$ with $k+l =n$.
Then for all integers $n > 81$, a simple observation of Holzer and K\"onig's proof of Theorem \ref{thm-max-prime} 
shows that $|U'| - |V| \geq \binom{n}{2}$.  For all integers $7\leq n\leq 81$, 
the fact that $|U'| - |V| \geq \binom{n}{2}$ has been verified computationally. It follows
from Corollary \ref{cor-root-monoid-bound} that if $M$ is not of the form $U_{k,l}$, we have
	\[\stcmp(\rt(L)) \leq |M| \leq |V| \leq |U_{k,l}| - \binom{n}{2}.\]

The fact that the bound is tight is an immediate consequence of Theorem \ref{thm-min-dfa}.
\end{proof}

If Conjecture \ref{con-ukl-max} is true, then for all $n\geq7$, where $n$ is not prime, the construction 
of ${\mathcal A}_{\Sigma,X_{k,l}}$ yields a language that is within $\binom{n}{2}$ of the maximum blow-up.  
We conjecture that this construction achieves the maximum.

\begin{conjecture}
For and integer $n\geq7$, there exist positive, coprime integers $k\geq2$, $l\geq3$, with $k+l=n$, 
such that if $L$ is a language over an alphabet of size 2, with $\stcmp(L) \leq n$, then
$\stcmp(\rt(L)) \leq |U_{k,l}| - \binom{n}{2}$. This bound is tight.
\end{conjecture}

The results concerning the largest monoid on $\geq 3$ generators are definite and much simpler. 
For this reason, on alphabets of size $\geq3$ we are able to give a much better bound.\\

\begin{lemma}
\label{lem-monoid-gap}
For $n \geq 1$, if $M \subseteq T_n$ is a monoid such that $|M| > n^n - \binom{n}{2}$,
then $M = T_n$.
\end{lemma}

\begin{proof}
For $1\leq n \leq 3$, the result can easily be verified computationally, so assume
that $n \geq 4$.\\

There are $\binom{n}{2}$ transpositions in $T_n$. Since $|M| > |T_n| - \binom{n}{2}$, it follows
that $M$ contains at least one transposition.  There are $(n-1)!$ permutations of $Z_n$ that have 
one cycle of length $n$.  Since $n \geq 4$, we have that $(n-1)! \geq \binom{n}{2}$.  Again,
considering the size of $M$, it follows that $M$ contains at least one  permutation 
that is a full $n$-cycle. It follows that $S_n \subseteq M$.\\

Furthermore, there are $\binom{n}{2}\cdot n!$ transformations of $Z_n$ that have rank $n-1$, so it follows 
that $M$ contains at least one transformation of rank $n-1$.  Then by Lemma \ref{lem-semi-gen}, 
we have that $M = T_n$.
\end{proof}

We now state our third main result.

\begin{theorem}
\label{thm-3-max}
Let $\Sigma$ be an alphabet of size $m \geq 3$. For $n \geq 1$, if $L$ is a language over $\Sigma$ 
with $\stcmp(L) \leq n$, then $\stcmp(\rt(L)) \leq n^n - \binom{n}{2}$.  Furthermore, 
this bound is tight.
\end{theorem}

\begin{proof}
Define $M$ to be the transformation monoid of the smallest DFA recognizing $L$. If
$|M| \leq n^n - \binom{n}{2}$, then certainly  $\stcmp(\rt(L)) \leq n^n - \binom{n}{2}$. So
suppose that $|M| > n^n - \binom{n}{2}$. Then it follows from Lemma \ref{lem-monoid-gap}
that $M = T_n$.\\

For $1 \leq n \leq 6$, it has been verified computationally that if the transformation monoid
of the minimal DFA recognizing $L$ is $T_n$, then $\stcmp(\rt(L)) = n^n - \binom{n}{2}$.
For $n \geq 7$, if the transformation monoid is $T_n$, then clearly $U_{k,l} \subseteq T_n$ for 
some suitable $k,l$ so that Theorem \ref{thm-min-dfa} applies, and hence $\stcmp(\rt(L)) = n^n - \binom{n}{2}$.\\

To show that the bound is tight, it suffices to show that for any $n$ there exists a language $L$ 
over $\Sigma$ such that the transformation monoid of the minimal DFA recognizing $L$ is $T_n$.
Let $X$ be a set of transformations such that $|X| = \min(n,3)$ and $X$ generates $T_n$. 
For $n\in\{1,2\}$, the fact that such an $X$ exists is easy to check.  For
$n\geq3$, the existence of $X$ follows from Lemma \ref{lem-semi-gen}. Then the language 
$L({\mathcal A}_{\Sigma,X})$, gives the desired result.
\end{proof}


\newpage
\appendix
{\Large\bf Appendix: Omitted Proofs}\\
\bigskip

\begin{proof}[Proof (Lemma \ref{lem-max-by-gap})]
As stated in \cite{HK03}, for $k + l = n$, we have the following formula
\begin{align*}
	|U_{k,l}| &= kl + \sum_{i=1}^{n}\biggl(\binom{n}{i}-\binom{k}{i-l}\biggr)
		\biggl(\stirlingtwo{n}{i}-\sum_{r=1}^{i}\stirlingtwo{k}{r}\stirlingtwo{l}{i-r}\biggr) i!,
\end{align*}
where $\stirlingtwo{n}{i}$ is a Stirling number of the second kind, the number
of ways to partition a set of $n$ elements into $i$ non-empty sets. This gives
\begin{align*}
	|U_{k,l}| - |U_{l,k}| &= 
		\sum_{i=1}^{n}\biggl(\binom{l}{i-k}-\binom{k}{i-l}\biggr)
		\biggl(\stirlingtwo{n}{i}-\sum_{r=1}^{i}\stirlingtwo{k}{r}\stirlingtwo{l}{i-r}\biggr)i!. \tag{$*$}
\end{align*}
Since $\stirlingtwo{n}{k} = 0$ whenever $k>n$ or $k<1$, for $k = 2$, and $l = n-2$, we have
\begin{align*}
	\sum_{r=1}^{i}\stirlingtwo{2}{r}\stirlingtwo{n-2}{i-r} &= \stirlingtwo{n-2}{i-2}+\stirlingtwo{n-2}{i-1}.
\end{align*}
Also, notice that $\binom{n-2}{i-2}-\binom{2}{i-n+2}$ is positive when $2\geq i\geq n-1$, and
zero otherwise, so that $(*)$ becomes
\begin{align*}
	|U_{2,n-2}| - |U_{n-2,2}| &\geq \sum_{i=2}^{n-1}
		\biggl(\stirlingtwo{n}{i}-\stirlingtwo{n-2}{i-2}-\stirlingtwo{n-2}{i-1}\biggr)i!.
\end{align*}
And finally, using the identity $\stirlingtwo{n}{k} = \stirlingtwo{n-1}{k-1} + k\stirlingtwo{n-1}{k}$, we
see that
\begin{align*}
	\stirlingtwo{n}{i} &= \stirlingtwo{n-2}{i-2} + (2i -1)\stirlingtwo{n-2}{i-1} + (i-1)\stirlingtwo{n-2}{i},
\intertext{so that we get}
	|U_{2,n-2}| - |U_{n-2,2}| &\geq \sum_{i=2}^{n-1}
		\biggl((2i-2)\stirlingtwo{n-2}{i-1}+(i-1)\stirlingtwo{n-2}{i}\biggr)i!
		\geq \sum_{i=2}^{n-1}i! \geq \binom{n}{2}.
\end{align*}
\end{proof}

\end{document}